\documentclass[english]{article}
\usepackage[T1]{fontenc}
\usepackage{textcomp}
\usepackage[latin9]{inputenc}
\usepackage{geometry}
\usepackage{amsmath}
\usepackage{amssymb}
\usepackage{esint}

\makeatletter
\newcommand{\lyxaddress}[1]{
	\par {\raggedright #1
	\vspace{1.4em}
	\noindent\par}
}

\makeatother

\makeatother

\usepackage{babel}

\makeatother

\usepackage{babel}

\makeatother

\usepackage{babel}

\makeatother

\usepackage{babel}

\makeatother

\usepackage{babel}
\begin{document}
\title{$\,_{3}F_{4}$ hypergeometric functions as a sum of a product of $\,_{2}F_{3}$
functions}
\author{Jack C. Straton}

\maketitle

\lyxaddress{Department of Physics, Portland State University, Portland, OR 97207-0751,
USA}

\lyxaddress{straton@pdx.edu}
\begin{abstract}
This paper shows that certain $\,_{3}F_{4}$ hypergeometric functions
can be expanded in sums of pair products of $\,_{1}F_{2}$ functions.
In special cases, the $\,_{3}F_{4}$ hypergeometric functions reduce
to $\,_{2}F_{3}$ functions. Further special cases allow one to reduce
the $\,_{2}F_{3}$ functions to $\,_{1}F_{2}$ functions, and the
sums to products of $\,_{0}F_{1}$ (Bessel) and $\,_{1}F_{2}$ functions.
This expands the class of hypergeometric functions having summation
theorems beyond those expressible as pair-products of generalized
Whittaker functions, $\,_{2}F_{1}$ functions, and $\,_{3}F_{2}$
functions into the realm of $\,_{p}F_{q}$ functions where $p<q$
for both the summand and terms in the series. 
\end{abstract}
\vspace{2pc}
 \\

\section{Introduction}

Jackson~\cite{Jackson,JacksonII} extended, in 1942, the concept of
a \emph{summation theorem} over pair-products of functions such as
Bessel functions~\cite{GR5}
( p. 992 No. 8.53) to sums of pair-products of generalized hypergeometric
functions, a much the broader class of functions. His main focus was
on hypergeometric functions of two variables (x, y), but as a special
case he replaced $y=xq^{b}$, which gave $\,_{2}F_{1}$ functions
expanded as pairs of $\,_{2}F_{1}$ functions (his Equation~(I.55))
and $\,_{1}F_{1}$ functions, also known as Whittaker functions, expanded
as pairs of $\,_{1}F_{1}$ functions (his Equation~(II.69)). Ragab
\cite{Ragab} found six such expressions in 1962 that involved Slater's
\cite{Slater} generalization of Whittaker functions to $\,_{p}F_{p}$
functions for larger values of \emph{p}. In all but one case, the
left-hand side has $x^{2}$ rather than $x$ as the argument in the
sum, such~as

\noindent
\begin{eqnarray}
\,_{2}F_{3} &  & \hspace{-1cm}\left(\frac{b}{2}+\frac{1}{2},\frac{b}{2};a+\frac{1}{2},\frac{c}{2}+\frac{1}{2},\frac{c}{2};\frac{x^{2}}{16}\right)=\sum_{r=0}^{\infty}\frac{(-1)^{r}x^{2r}(a)_{r}(b)_{r}(c-b)_{r}}{r!(2a)_{r}(c)_{2r}(c+r-1)_{r}}\nonumber \\
 & \times & \,_{2}F_{2}(a+r,b+r;2a+r,c+2r;x)\,_{1}F_{1}\left(b+r;c+2r;-\frac{x}{2}\right)\:.\label{eq:Ragab}
\end{eqnarray}

In addition to such sums, he also expressed a $\,_{4}F_{5}$ function
as a sum of products of $\,_{2}F_{2}$ functions. In 1964, Verma \cite{Verma}
rederived Jackson's and some of Ragab's results and added expansions
of $\,_{3}F_{2}$ functions as a product of a $\,_{2}F_{1}$ function
with another $\,_{3}F_{2}$ function. He also expressed a $\,_{5}F_{5}$
generalized Whittaker function as a sum of products of $\,_{2}F_{2}$
functions.

There is also an active modern field of study of summations theorems
not involving \emph{pair-products} of generalized hypergeometric functions
in the summands. To name a few of many, Choi, Milovanovi and Rathie~\cite{Choi}
express Kamp{\'e} de F{\'e}riet functions as finite sums, as do Wang and Chen~\cite{WangChen}.
Wang~\cite{Wang} gives Kamp{\'e} de F{\'e}riet functions and various related
functions as infinite sums. Yakubovich \cite{Yakubovich} does likewise
for generalized hypergeometric functions. Liu and Wang \cite{LiuWang}
were able to reduce Kamp{\'e} de F{\'e}riet functions to Appell series and
generalized hypergeometric functions by generalizing of classical
summation theorems that have been done by Kummer, Gauss, and Bailey,
work extended by Choi and Rathie. \cite{ChoiRathie} For an excellent
summary, and extensions, of those classical summation theorems for
generalized hypergeometric functions, please see Awad et al~\cite{Awad} 

In a prior paper, \cite{stra24c} the author returned to infinite
summation theorems for hypergeometric functions in terms of \emph{pair-products}
of hypergeometric functions, specifically $\,_{3}F_{4}$ hypergeometric
functions expressed as in terms of pair-products\emph{ }of $\,_{2}F_{3}$
hypergeometric functions. Special values of their parameters reduced
these to $\,_{2}F_{3}$ functions expanded in sums of pair products
of $\,_{1}F_{2}$ functions. While interesting in itself, this result
had a specific application in calculating the response of atoms to
laser stimulation in the Strong Field Approximation (SFA) \cite{Reiss1,Reiss3,Reiss4,Reiss5,Reiss6,Faisel2}.
Whereas perturbation expansions will not converge if the applied laser
field is sufficiently large, the Strong Field Approximation (SFA)
is an analytical approximation that is non-perturbative. Keating \cite{KeatingPhD}
applied it specifically to the production of the positive antihydrogen
ion.

The present work builds on that paper. Section 2 lays out the prior
approach in a clean pattern divested of the terminology of the SFA,
and the fact that it centers on the transition amplitude for the positive
antihydrogen ion. We simply evaluate the integral 

\begin{equation}
A_{\mu\nu}\left(k\right)=\int_{-1}^{1}J_{\mu}\left(k\,u\right)J_{\nu}\left(k\,u\right)\,du\label{eq:Legendre_JJ}
\end{equation}
two ways, first by expanding each Bessel function in a series of Legendre
polynomials,

\begin{equation}
J_{N}\left(kx\right)=\sum_{L=0}^{\infty}a_{LN}\left(k\right)P_{L}(x)\:\label{eq:Legendre=000020series}
\end{equation}
and comparing that to direct integration.

Section 3 steps entirely away from an integral with a known physical
application, by inserting complicating factors in (\ref{eq:Legendre_JJ})
such that it is easily integrated using series of Chebyshev polynomials
instead of Legendre ones. Section 4 inserts other complicating factors
in (\ref{eq:Legendre_JJ}) so that it is readily solved using series
of Gegenbauer polynomials. The result in each case is a set of $\,_{3}F_{4}$
hypergeometric functions expressed as in terms of pair-products\emph{
}of $\,_{1}F_{2}$ hypergeometric functions.

\section{Review of the Legendre version}

The coefficients of the Fourier-Legendre series for the Bessel function
$J_{N}\left(kx\right)$ (\ref{eq:Legendre=000020series}) were given
in Jet Wimp's 1962 Jacobi expansion \cite{jacobi_cheb_Wimp} of the
Anger-Weber function (his equations (2.10) and (2.11)) since Legendre
polynomials are a subset of Jacobi polynomials and the Bessel function
$J_{N}\left(kx\right)$ are a special case of the Anger-Weber function
$\mathbf{J}_{\nu}\left(kx\right)$ when $\nu$ is an integer, the
first two lines of

\begin{equation}
\begin{array}{ccc}
a_{LN}\left(k\right) & = & \frac{\sqrt{\pi}2^{-2L-2}(2L+1)k^{L}i^{L-N}}{\Gamma\left(\frac{1}{2}(2L+3)\right)}\left(1+(-1)^{L+N}\right)\binom{L}{\frac{L-N}{2}}\\
 & \times & \,_{2}F_{3}\left(\frac{L}{2}+\frac{1}{2},\frac{L}{2}+1;L+\frac{3}{2},\frac{L}{2}-\frac{N}{2}+1,\frac{L}{2}+\frac{N}{2}+1;-\frac{k^{2}}{4}\right)\\
 & = & \sqrt{\pi}\,2^{-2L-2}(2L+1)k^{L}i^{L-N}\left(1+(-1)^{L+N}\right)\Gamma(L+1)\\
 & \times & \,_{2}\tilde{F}_{3}\left(\frac{L}{2}+\frac{1}{2},\frac{L}{2}+1;L+\frac{3}{2},\frac{L}{2}-\frac{N}{2}+1,\frac{L}{2}+\frac{N}{2}+1;-\frac{k^{2}}{4}\right)\:.\\
 & = & \sqrt{\pi}(2L+1)2^{-L-1}i^{L-\text{N}}\sum_{M=0}^{\infty}\frac{\left(\left(-\frac{1}{4}\right)^{M}k^{L+2M}\right)}{2^{L+2M+1}\left(M!\Gamma\left(L+M+\frac{3}{2}\right)\right)}\\
 & \times & \left(1+(-1)^{L+2M+\text{N}}\right)\binom{L+2M}{\frac{1}{2}(L+2M-\text{N})}
\end{array}\label{eq:a_L_as_2F3}
\end{equation}
In numerical checks, I found that the $\,_{2}F_{3}$ hypergeometric
function becomes infinite whenever \emph{N>1} is an integer larger
than \emph{L}, while the binomial prefactor simultaneously goes to
zero. Wimp's work is entirely analytical and does not mention any
such calculational difficulties. They may be removed by introducing
regularized hypergeometric functions \cite{wolfram.com/07.26.26.0001.01}

\begin{equation}
\,_{2}F_{3}\left(a_{1},a_{2};b_{1},b_{2},b_{3};z\right)=\Gamma\left(b_{1}\right)\Gamma\left(b_{2}\right)\Gamma\left(b_{3}\right)\,_{2}\tilde{F}_{3}\left(a_{1},a_{2};b_{1},b_{2},b_{3};z\right)\label{eq:regularized_2F3}
\end{equation}
and cancelling the $\Gamma\left(b_{i}\right)$ with the binomial function
to remove such indeterminacies in the computation, giving the second
form in (\ref{eq:a_L_as_2F3}). The third form in (\ref{eq:a_L_as_2F3})
is from Keating. \cite{KeatingPhD}

With the expansion (\ref{eq:Legendre=000020series}), and using the
orthogonality of the Legendre polynomials to reduce the double sum
into a single sum over a product of $\,_{2}F_{3}$ hypergeometric
functions, one has
\begin{equation}
\begin{array}{ccc}
A_{\mu\nu}\left(k\right) & = & {\displaystyle \sum_{L=0}^{\infty}\frac{\pi^{2}2^{-8L-5}(2L+1)k^{2L}\left((-1)^{L+\mu}+1\right)\left((-1)^{L+\nu}+1\right)\Gamma(2L+2)^{2}i^{2L-\mu-\nu}}{\Gamma\left(\frac{1}{2}(2L+3)\right)^{2}}}\\
 & \times & \,_{2}\tilde{F}_{3}\left(\frac{L}{2}+\frac{1}{2},\frac{L}{2}+1;L+\frac{3}{2},\frac{L}{2}-\frac{\mu}{2}+1,\frac{L}{2}+\frac{\mu}{2}+1;-\frac{k^{2}}{4}\right)\\
 & \times & \,_{2}\tilde{F}_{3}\left(\frac{L}{2}+\frac{1}{2},\frac{L}{2}+1;L+\frac{3}{2},\frac{L}{2}-\frac{\nu}{2}+1,\frac{L}{2}+\frac{\nu}{2}+1;-\frac{k^{2}}{4}\right)\\
 & = & {\displaystyle \sum_{L=0}^{\infty}\frac{\pi^{2}2^{-8L-5}(2L+1)k^{2L}\left((-1)^{L+\mu}+1\right)\left((-1)^{L+\nu}+1\right)\Gamma(2L+2)^{2}i^{2L-\mu-\nu}}{\Gamma\left(\frac{1}{2}(2L+3)\right)^{2}}}\\
 &  & {\displaystyle \frac{1}{\Gamma\left(L+\frac{3}{2}\right)^{2}\Gamma\left(\frac{L}{2}-\frac{\mu}{2}+1\right)\Gamma\left(\frac{L}{2}+\frac{\mu}{2}+1\right)\Gamma\left(\frac{L}{2}-\frac{\nu}{2}+1\right)\Gamma\left(\frac{L}{2}+\frac{\nu}{2}+1\right)}}\\
 & \times & \,_{2}F_{3}\left(\frac{L}{2}+\frac{1}{2},\frac{L}{2}+1;L+\frac{3}{2},\frac{L}{2}-\frac{\mu}{2}+1,\frac{L}{2}+\frac{\mu}{2}+1;-\frac{k^{2}}{4}\right)\\
 & \times & \,_{2}F_{3}\left(\frac{L}{2}+\frac{1}{2},\frac{L}{2}+1;L+\frac{3}{2},\frac{L}{2}-\frac{\nu}{2}+1,\frac{L}{2}+\frac{\nu}{2}+1;-\frac{k^{2}}{4}\right)
\end{array}\:.\label{eq:Legendre_2F3_2F3}
\end{equation}

Alternatively, one can transform from Bessel to hypergeometric functions
using \cite{PBM3} (p. 220 No. 2.21.2.11),
\cite{Luke1} (p. 212 No.
6.2.7.1), or \cite{wolfram.com/03.01.26.0002.01} 
\begin{equation}
J_{\nu}(z)=\frac{\left(\frac{z}{2}\right)^{\nu}\,_{0}F_{1}\left(;\nu+1;-\frac{z^{2}}{4}\right)}{\Gamma(\nu+1)}\label{eq:J_to_0F1}
\end{equation}
 and combine pairs using \cite{EH,wolfram.com/07.17.26.0004.01}
or \cite{Luke1} (p. 228
No. 6.4.1.26) 
\begin{equation}
\,_{0}F_{1}(;b;z)\,_{0}F_{1}(;c;z)=\,_{2}F_{3}\left(\frac{b}{2}+\frac{c}{2}-\frac{1}{2},\frac{b}{2}+\frac{c}{2};b,c,b+c-1;4z\right)\label{eq:0F1=0000200F1=000020to=0000202F3}
\end{equation}
 so that \cite{Luke1}
(p. 216 No. 6.2.7.39) 
\begin{equation}
J_{\mu}(z)J_{\nu}(z)=\frac{2^{-\mu-\nu}z^{\mu+\nu}}{\Gamma(\mu+1)\Gamma(\nu+1)}\,_{2}F_{3}\left(\frac{\mu}{2}+\frac{\nu}{2}+\frac{1}{2},\frac{\mu}{2}+\frac{\nu}{2}+1;\mu+1,\nu+1,\mu+\nu+1;-z^{2}\right)\:.\label{eq:JJ=000020to=0000202F3}
\end{equation}

\noindent Since the integral we wish to perform has symmetrical limits,
only even integrands ($\mu+\nu$ even) will be nonzero: 
\begin{eqnarray}
A_{\mu\nu}\left(k\right) & = & \int_{-1}^{1}J_{\mu}\left(k\,u\right)J_{\nu}\left(k\,u\right)\,du\nonumber \\
 & = & \int_{-1}^{1}\frac{2^{-\mu-\nu}(ku)^{\mu+\nu}}{\Gamma(\mu+1)\Gamma(\nu+1)}\,_{2}F_{3}\left(\frac{\mu}{2}+\frac{\nu}{2}+\frac{1}{2},\frac{\mu}{2}+\frac{\nu}{2}+1;\mu+1,\nu+1,\mu+\nu+1;-k^{2}u^{2}\right)\,du\label{eq:convert=000020to=000020y}\\
 & \equiv & \int_{-1}^{1}f\left(u^{2}\right)du=2\int_{0}^{1}f\left(u^{2}\right)du=\int_{0}^{1}f\left(y\right)y^{-1/2}dy\:.\nonumber 
\end{eqnarray}
 We can then use \cite{PBM3} (p.
334 No. 2.22.2.1) 
\begin{align}
\int_{0}^{a}y^{\alpha-1} & (a-y)^{\beta-1}\,_{p}F_{q}(a_{1},\ldots,\,a_{p},b_{1},\ldots,\,b_{q};-\omega y)\,dy\label{eq:pFq_integral}\\
 & =\frac{\Gamma(\alpha)\Gamma(\beta)a^{\alpha+\beta-1}}{\Gamma(\alpha+\beta)}\,_{p+1}F_{q+1}(a_{1},\ldots,\,a_{p},\alpha;b_{1},\ldots,\,b_{q},\alpha+\beta;-a\omega)\nonumber \\
 & \hspace{2.6cm}\left[\Re(\alpha)>0\land\Re(\beta)>0\land a>0\right]\nonumber 
\end{align}
 with $\alpha=\mu+\nu$ even, $a=1$, and $\beta=1$ to find

\noindent
\begin{eqnarray}
A_{\mu\nu}\left(k\right) & = & \frac{2^{-\mu-\nu}k^{\mu+\nu}}{(\mu+\nu+1)\Gamma(\mu+1)\Gamma(\nu+1)}\nonumber \\
 & \times & \,_{3}F_{4}\left(\frac{\mu}{2}+\frac{\nu}{2}+\frac{1}{2},\frac{\mu}{2}+\frac{\nu}{2}+\frac{1}{2},\frac{\mu}{2}+\frac{\nu}{2}+1;\mu+1,\frac{\mu}{2}+\frac{\nu}{2}+\frac{3}{2},\nu+1,\mu+\nu+1;-k^{2}\right)\label{eq:Legendre_3F4}
\end{eqnarray}
which is the desired summation of (\ref{eq:Legendre_2F3_2F3}). One
may formalize this relation:

\noindent\textbf{Theorem 1.} For integer $\mu$, and $\nu$, $\mu+\nu$
even, but for any values of \emph{k}, 

\begin{equation}
\begin{array}{ccc}
 &  & \,_{3}F_{4}\left(\frac{\mu}{2}+\frac{\nu}{2}+\frac{1}{2},\frac{\mu}{2}+\frac{\nu}{2}+\frac{1}{2},\frac{\mu}{2}+\frac{\nu}{2}+1;\mu+1,\frac{\mu}{2}+\frac{\nu}{2}+\frac{3}{2},\nu+1,\mu+\nu+1;-k^{2}\right)\\
 &  & \\
 & = & 2^{\mu+\nu-1}(\mu+\nu+1)\Gamma(\mu+1)\Gamma(\nu+1)k^{-\mu-\nu}\\
 & \times & {\displaystyle \sum_{L=0}^{\infty}\frac{\pi^{2}2^{-8L-5}(2L+1)k^{2L}\left((-1)^{L+\mu}+1\right)\left((-1)^{L+\nu}+1\right)\Gamma(2L+2)^{2}i^{2L-\mu-\nu}}{\Gamma\left(\frac{1}{2}(2L+3)\right)^{2}}}\\
 & \times & \,_{2}\tilde{F}_{3}\left(\frac{L}{2}+\frac{1}{2},\frac{L}{2}+1;L+\frac{3}{2},\frac{L}{2}-\frac{\mu}{2}+1,\frac{L}{2}+\frac{\mu}{2}+1;-\frac{k^{2}}{4}\right)\:.\\
 & \times & \hspace{-0cm}\,_{2}\tilde{F}_{3}\left(\frac{L}{2}+\frac{1}{2},\frac{L}{2}+1;L+\frac{3}{2},\frac{L}{2}-\frac{\nu}{2}+1,\frac{L}{2}+\frac{\nu}{2}+1;-\frac{k^{2}}{4}\right)\\
 &  &\\
 & = & 2^{\mu+\nu-1}(\mu+\nu+1)\Gamma(\mu+1)\Gamma(\nu+1)k^{-\mu-\nu}\\
 & \times & {\displaystyle \sum_{L=0}^{\infty}\frac{\pi^{2}2^{-8L-5}(2L+1)k^{2L}\left((-1)^{L+\mu}+1\right)\left((-1)^{L+\nu}+1\right)\Gamma(2L+2)^{2}i^{2L-\mu-\nu}}{\Gamma\left(\frac{1}{2}(2L+3)\right)^{2}}}\\
 &  & {\displaystyle \frac{1}{\Gamma\left(L+\frac{3}{2}\right)^{2}\Gamma\left(\frac{L}{2}-\frac{\mu}{2}+1\right)\Gamma\left(\frac{L}{2}+\frac{\mu}{2}+1\right)\Gamma\left(\frac{L}{2}-\frac{\nu}{2}+1\right)\Gamma\left(\frac{L}{2}+\frac{\nu}{2}+1\right)}}\\
 & \times & \,_{2}F_{3}\left(\frac{L}{2}+\frac{1}{2},\frac{L}{2}+1;L+\frac{3}{2},\frac{L}{2}-\frac{\mu}{2}+1,\frac{L}{2}+\frac{\mu}{2}+1;-\frac{k^{2}}{4}\right)\\
 & \times & \hspace{-0cm}\,_{2}F_{3}\left(\frac{L}{2}+\frac{1}{2},\frac{L}{2}+1;L+\frac{3}{2},\frac{L}{2}-\frac{\nu}{2}+1,\frac{L}{2}+\frac{\nu}{2}+1;-\frac{k^{2}}{4}\right)
\end{array}\label{eq:theorem_1}
\end{equation}

\section{Chebyshev version}

We turn now to the proof of a theorem based on series of Chebyshev
polynomials:

\noindent\textbf{Theorem 2.} For any values of \emph{k}, $\mu$,
and $\nu$, 

\begin{equation}
\begin{array}{ccc}
 &  & 
 \,_{3}F_{4}\left(\frac{1}{2},\frac{\mu}{2}+\frac{\nu}{2}+\frac{1}{2},\frac{\mu}{2}+\frac{\nu}{2}+1;1,\mu+1,\nu+1,\mu+\nu+1;-k^{2}\right)\\
 &  & \\
 & = & 2^{\mu+\nu}\Gamma(\mu+1)\Gamma(\nu+1){\displaystyle \sum_{L=0}^{\infty}\frac{(-1)^{2L}k^{4L}2^{-8L-\mu-\nu}\left(2-\delta_{L0}\right)}{(L!)^{2}\Gamma(L+\mu+1)\Gamma(L+\nu+1)}}\\
 & \times & \,_{1}F_{2}\left(L+\frac{1}{2};2L+1,L+\mu+1;-\frac{k^{2}}{4}\right)\:.\\
 & \times & \hspace{-0cm}\,_{1}F_{2}\left(L+\frac{1}{2};2L+1,L+\nu+1;-\frac{k^{2}}{4}\right)
\end{array}\label{eq:theorem_2}
\end{equation}

\textbf{Proof of Theorem 2.} 

Wimp \cite{jacobi_cheb_Wimp} gives an expansion of the more complicated function $J_{\nu}\left(kx\right)\left(kx\right)^{-\nu}$
in a series of Chebyshev polynomials 

\begin{equation}
J_{\nu}\left(kx\right)\left(kx\right)^{-\nu}=\sum_{L=0}^{\infty}C_{L\nu}\left(k\right)T_{2L}(x)\:\quad-1\leq x\leq1\label{eq:Chebyshev=000020series}
\end{equation}
(which applies to non-integer indices as well), where 

\begin{equation}
\begin{array}{ccc}
C_{L\nu}\left(k\right) & = & {\displaystyle \frac{(-1)^{L}k^{2L}2^{-4L-\nu}\left(2-\delta_{L0}\right)}{L!\Gamma(L+\nu+1)}}\,_{1}F_{2}\left(L+\frac{1}{2};L+\nu+1,2L+1;-\frac{k^{2}}{4}\right)\end{array}\;.\label{eq:C_L_as_1F2}
\end{equation}
So let us consider an integral akin to (\ref{eq:Legendre_JJ}), 

\begin{equation}
B_{\mu\nu}\left(k\right)=\int_{-1}^{1}J_{\mu}\left(k\,u\right)\left(ku\right)^{-\mu}J_{\nu}\left(k\,u\right)\left(ku\right)^{-\nu}\frac{1}{\sqrt{1-u^{2}}}\,du\:,\label{eq:Chebyshev_JJ}
\end{equation}
that is complicated not only by these inverse powers, but by the weight
function that appears in the orthogonality relation for Chebyshev
polynomials \cite{GR5}
(p. 1057 No. 8.949.9):

\begin{equation}
\int_{-1}^{1}T_{n}(x)T_{m}(x)\frac{1}{\sqrt{1-x^{2}}}=\begin{cases}
\begin{array}{c}
0\\
\pi/2\\
\pi
\end{array} & \begin{array}{c}
m\neq n\\
m=n\neq0\\
m=n=0
\end{array}\end{cases}\:.\label{eq:T_orth}
\end{equation}

Again we evaluate (\ref{eq:Chebyshev_JJ}) two ways. One may first
expand each factor (\ref{eq:Chebyshev=000020series}) in (\ref{eq:Chebyshev_JJ})
and use the orthogonality of the Chebyshev polynomials, above, to
reduce the double sum into a single sum over a product of $\,_{1}F_{2}$
hypergeometric functions.

\begin{equation}
\begin{array}{ccc}
B_{\mu\nu}\left(k\right) & = & \pi{\displaystyle \sum_{L=0}^{\infty}\frac{(-1)^{2L}k^{4L}2^{-8L-\mu-\nu}\left(2-\delta_{L0}\right)}{(L!)^{2}\Gamma(L+\mu+1)\Gamma(L+\nu+1)}}\\
 & \times & \,_{1}F_{2}\left(L+\frac{1}{2};2L+1,L+\mu+1;-\frac{k^{2}}{4}\right)\,_{1}F_{2}\left(L+\frac{1}{2};2L+1,L+\nu+1;-\frac{k^{2}}{4}\right)
\end{array}\label{eq:B_as_1F2_1F2}
\end{equation}

One may instead convert the Bessel function product to a $\,_{2}F_{3}$
hypergeometric function using (\ref{eq:JJ=000020to=0000202F3}). In
the present case, the leftover powers of $\left(ku\right)^{-\mu-\nu}$
in (\ref{eq:Chebyshev_JJ}) precisely cancel the powers in (\ref{eq:JJ=000020to=0000202F3}),
so that the only power remaining after we convert from \emph{u} to
\emph{y} in (\ref{eq:convert=000020to=000020y}) is the square root
of \emph{y} in the denominator so that $\alpha=\frac{1}{2}$ . Again
$a=1$, but we now have the factor $\sqrt{1-y}$ in the denominator
so that$\beta=\frac{1}{2}$. Thus,

\noindent
\begin{eqnarray}
B_{\mu\nu}\left(k\right) & = & \frac{\pi2^{-\mu-\nu}}{\Gamma(\mu+1)\Gamma(\nu+1)}\,_{3}F_{4}\left(\frac{1}{2},\frac{\mu}{2}+\frac{\nu}{2}+\frac{1}{2},\frac{\mu}{2}+\frac{\nu}{2}+1;1,\mu+1,\nu+1,\mu+\nu+1;-k^{2}\right)\label{eq:Chebyshev_3F4}
\end{eqnarray}
which completes the proof. $\square$

While this derivation formally holds for any values of $\mu$ and
$\nu$, those who wish to utilize this in numerical calculations will
find that if one or both of $\mu$ or $\nu$ is a negative integer
of sufficient magnitude, one will likely need to introduce 
 regularized hypergeometric functions to cancel the gamma functions
in the denominators to avoid infinities and/or indeteminancies. 

When $\nu=-\mu$ the the order of the hypergeometric function on the
left-hand side of (\ref{eq:theorem_2}) is reduced since one upper
parameter equals a lower one (equaling one). This will transform
(\ref{eq:Chebyshev_3F4}) to

\noindent
\begin{eqnarray}
\hspace{-0.6 cm} B_{\mu,-\mu}\left(k\right) & = & \frac{\pi}{\Gamma(1-\mu)\Gamma(\mu+1)}\,_{2}F_{3}\left(\frac{1}{2},\frac{1}{2};1,1-\mu,\mu+1;-k^{2}\right)=\pi\,_{2}\tilde{F}_{3}\left(\frac{1}{2},\frac{1}{2};1,1-\mu,\mu+1;-k^{2}\right),\label{eq:Chebyshev_2F3}
\end{eqnarray}
where the second form involving the regularized hypergeometric function
\cite{wolfram.com/07.26.26.0001.01} will be needed in numerical evaluations
to avoid infinities and/or indeteminancies whenever $\mu$ is a nonzero
integer. The corresponding gamma functions in the denominator of (\ref{eq:B_as_1F2_1F2})
will likewise cause numerical problems unless we switch to regularized
hypergeometric functions,\cite{wolfram.com/07.22.26.0001.01}

\begin{equation}
\,_{1}F_{2}\left(a_{1};b_{1},b_{2};z\right)=\Gamma\left(b_{1}\right)\Gamma\left(b_{2}\right)\,_{1}\tilde{F}_{2}\left(a_{1};b_{1},b_{2};z\right)\:.\label{eq:regularized_1F2}
\end{equation}
These give

\noindent\textbf{Corollary 2a.} For any values of \emph{k} and $\mu$ 

\begin{equation}
\begin{array}{ccc}
\,_{2}\tilde{F}_{3}\left(\frac{1}{2},\frac{1}{2};1,1-\mu,\mu+1;-k^{2}\right) & = & {\displaystyle \sum_{L=0}^{\infty}\frac{(-1)^{2L}k^{4L}\Gamma(2L+1)^{2}2^{-8L-\mu-\nu}\left(2-\delta_{L0}\right)}{(L!)^{2}}}\\
 & \times & \,_{1}\tilde{F}_{2}\left(L+\frac{1}{2};2L+1,L+\mu+1;-\frac{k^{2}}{4}\right)\:.\\
 & \times & \,_{1}\tilde{F}_{2}\left(L+\frac{1}{2};2L+1,L+\nu+1;-\frac{k^{2}}{4}\right)
\end{array}\label{eq:Corollary_2a}
\end{equation}

When either of $\mu$ or $\nu$ is $-\frac{1}{2}$ the order of the
hypergeometric function on the left-hand side of (\ref{eq:theorem_2})
is also reduced since the second or third lower parameter equals the
first upper parameter $\left(\frac{1}{2}\right)$:

\noindent\textbf{Corollary 2b.} For any values of \emph{k} and $\mu$ 

\begin{equation}
\begin{array}{ccc}
\,_{2}F_{3}\left(\frac{\mu}{2}+\frac{1}{4},\frac{\mu}{2}+\frac{3}{4};1,\mu+\frac{1}{2},\mu+1;-k^{2}\right) & = & \sqrt{\pi}\Gamma(\mu+1){\displaystyle \sum_{L=0}^{\infty}\frac{(-1)^{2L}2^{-6L}k^{2L}\Gamma(2L+1)\left(2-\delta_{L0}\right)}{(L!)^{2}\Gamma\left(L+\frac{1}{2}\right)\Gamma(L+\mu+1)}}\\
 & \times & J_{2L}(k)\hspace{-0cm}\,_{1}F_{2}\left(L+\frac{1}{2};2L+1,L+\mu+1;-\frac{k^{2}}{4}\right)\:.
\end{array}\label{eq:Corollary_2b}
\end{equation}
An alternative form may be had by replacing the Bessel function in
the second line by a $\,_{0}F_{1}$ hypergeometric function using (\ref{eq:J_to_0F1}).

\section{Gegenbauer version}

A similar theorem arises from integrating over a product of series
of Gegenbauer polynomials:

\noindent\textbf{Theorem 3.} For any values of \emph{k}, $\mu$,
$\nu$, and for $\lambda\neq0$,

\begin{equation}
\begin{array}{ccc}
&  &
\,_{3}F_{4}\left(\frac{1}{2},\frac{\mu}{2}+\frac{\nu}{2}+\frac{1}{2},\frac{\mu}{2}+\frac{\nu}{2}+1;\lambda+1,\mu+1,\nu+1,\mu+\nu+1;-k^{2}\right)\\
 &  &\\
 & = & {\displaystyle \frac{\Gamma(\lambda+1)2^{\mu+\nu}\Gamma(\mu+1)\Gamma(\nu+1)}{\sqrt{\pi}\Gamma\left(\lambda+\frac{1}{2}\right)\Gamma(\lambda)^{2}}}\\
 &  & {\displaystyle \sum_{L=0}^{\infty}\frac{(-1)^{2L}k^{4L}\left(\left(\lambda+\frac{1}{2}\right)_{2L}\right){}^{2}\Gamma(2L+2\lambda)2^{4L-2\lambda-\mu-\nu+1}}{(2L)!(2L+\lambda)\left((2\lambda)_{2L}\right){}^{2}\left((2L+2\lambda)_{2L}\right){}^{2}\left(L+\frac{1}{2}\right)_{\mu+\frac{1}{2}}\left(L+\frac{1}{2}\right)_{\nu+\frac{1}{2}}}}\\
 & \times & \,_{1}F_{2}\left(L+\frac{1}{2};2L+\lambda+1,L+\mu+1;-\frac{k^{2}}{4}\right)\:.\\
 & \times & \,_{1}F_{2}\left(L+\frac{1}{2};2L+\lambda+1,L+\nu+1;-\frac{k^{2}}{4}\right)
\end{array}\label{eq:theorem_3}
\end{equation}

\textbf{Proof of Theorem 3} 

Though Wimp  \cite{jacobi_cheb_Wimp} does not do so, in another prior work \cite{stra24d}  we  showed that one may use Wimp's Jacobi expansion to find Gegenbauer polynomial expansions of Bessel functions, since
\cite{GR5} (p. 1060 No.
8.962.4)
\begin{equation}
P_{2n}^{\left(\lambda-\frac{1}{2},\lambda-\frac{1}{2}\right)}(z)=\frac{\left(\lambda+\frac{1}{2}\right)_{2n}C_{2n}^{\lambda}(z)}{(2\lambda)_{2n}}\;.\label{eq:Jacobi_Gegenbauer}
\end{equation}
This leads one to the following expansion:

\begin{equation}
J_{\nu}\left(kx\right)\left(kx\right)^{-\nu}=\sum_{L=0}^{\infty}b_{L\nu}\left(k\right)C_{2L}^{\lambda}(x)\:,\quad-1\leq x\leq1\label{eq:Gegenbauer=000020series}
\end{equation}
(which applies to non-integer indices as well) where where the coefficients
are 

\begin{equation}
\begin{array}{ccc}
b_{L\nu}\left(k\right) & ={\displaystyle \frac{(-1)^{L}k^{2L}2^{2L-\nu}\left(\lambda+\frac{1}{2}\right)_{2L}}{\sqrt{\pi}(2\lambda)_{2L}(2L+2\lambda)_{2L}\left(L+\frac{1}{2}\right)_{\nu+\frac{1}{2}}}} & \,_{1}F_{2}\left(L+\frac{1}{2};2L+\lambda+1,L+\nu+1;-\frac{k^{2}}{4}\right)\end{array}\;.\label{eq:b_L_as_1F2=000020Geg}
\end{equation}

So let us consider an integral akin to (\ref{eq:Chebyshev_JJ}), 

\begin{equation}
H_{\mu\nu}^{\lambda}\left(k\right)=\int_{-1}^{1}J_{\mu}\left(k\,u\right)\left(ku\right)^{-\mu}J_{\nu}\left(k\,u\right)\left(ku\right)^{-\nu}\left(1-u^{2}\right)^{\lambda-\frac{1}{2}}\,du\:,\label{eq:Gegenbauer_JJ}
\end{equation}
but containing the weight function that appears in the orthogonality
relation for Gegenbauer polynomials \cite{GR5}
(p. 1054 No. 8.939.8):

\begin{equation}
\int_{-1}^{1}C_{n}^{\lambda}(x)C_{m}^{\lambda}(x)\left(1-x^{2}\right)^{\lambda-\frac{1}{2}}=\begin{cases}
\begin{array}{c}
0\\
{\displaystyle \frac{\pi2^{1-2\lambda}\Gamma(n+2\lambda)}{n!(n+\lambda)\Gamma(\lambda)^{2}}}
\end{array} & \begin{array}{c}
m\neq n\\
m=n
\end{array}\end{cases}\left(\lambda\neq0\right)\:.\label{eq:C_orth}
\end{equation}

Again we evaluate it two ways. One may first expand each factor (\ref{eq:Gegenbauer=000020series})
in (\ref{eq:Gegenbauer_JJ}) and use the orthogonality of the Gegenbauer
polynomials, above, to reduce the double sum into a single sum over
a product of $\,_{1}F_{2}$ hypergeometric functions:

\begin{equation}
\begin{array}{ccc}
H_{\mu\nu}^{\lambda}\left(k\right) & = & {\displaystyle \sum_{L=0}^{\infty}\frac{(-1)^{2L}k^{4L}\left(\left(\lambda+\frac{1}{2}\right)_{2L}\right){}^{2}\Gamma(2L+2\lambda)2^{4L-2\lambda-\mu-\nu+1}}{(2L)!(2L+\lambda)\Gamma(\lambda)^{2}\left((2\lambda)_{2L}\right){}^{2}\left((2L+2\lambda)_{2L}\right){}^{2}\left(L+\frac{1}{2}\right)_{\mu+\frac{1}{2}}\left(L+\frac{1}{2}\right)_{\nu+\frac{1}{2}}}}\\
 & \times & \,_{1}F_{2}\left(L+\frac{1}{2};2L+\lambda+1,L+\mu+1;-\frac{k^{2}}{4}\right)\,_{1}F_{2}\left(L+\frac{1}{2};2L+\lambda+1,L+\nu+1;-\frac{k^{2}}{4}\right) \;.
\end{array}\label{eq:H_as_1F2_1F2_Gegen}
\end{equation}

One may instead convert the Bessel function product to a $\,_{2}F_{3}$
hypergeometric function using (\ref{eq:JJ=000020to=0000202F3}). In
the present case, the leftover powers of $\left(ku\right)^{-\mu-\nu}$in
(\ref{eq:Gegenbauer_JJ}) precisely cancel the powers in (\ref{eq:JJ=000020to=0000202F3}),
so that the only power remaining after we convert from \emph{u} to
\emph{y} in (\ref{eq:convert=000020to=000020y}) is the square root
of \emph{y} in the denominator so that $\alpha=\frac{1}{2}$ . Again
$a=1$, but we now have the factor $(1-y)^{\lambda-\frac{1}{2}}$
in the numerator so that $\beta=\lambda+\frac{1}{2}$. Thus,

\noindent
\begin{eqnarray}
 \hspace{-0.6 cm} H_{\mu\nu}^{\lambda}\left(k\right)  \hspace{-0.3 cm} & = & \hspace{-0.3 cm}  \frac{\sqrt{\pi}\Gamma\left(\lambda+\frac{1}{2}\right)2^{-\mu-\nu}}{\Gamma(\lambda+1)\Gamma(\mu+1)\Gamma(\nu+1)}\,_{3}F_{4}\left(\frac{1}{2},\frac{\mu}{2}+\frac{\nu}{2}+\frac{1}{2},\frac{\mu}{2}+\frac{\nu}{2}+1;\lambda+1,\mu+1,\nu+1,\mu+\nu+1;-k^{2}\right) \hspace{-0.1 cm} ,\label{eq:Gegenbauer_3F4}
\end{eqnarray}
which completes the proof. $\square$

When $\nu=-\mu$ the the order of the hypergeometric function on the
left-hand side of (\ref{eq:theorem_3}) is reduced since one upper
parameter equals a lower one (equaling one). This will transform (\ref{eq:Gegenbauer_3F4})
to

\noindent
\begin{eqnarray}
H_{\mu\nu}^{\lambda}\left(k\right) & = & \frac{\sqrt{\pi}\Gamma\left(\lambda+\frac{1}{2}\right)}{\Gamma(\lambda+1)\Gamma(1-\mu)\Gamma(\mu+1)}\,_{2}F_{3}\left(\frac{1}{2},\frac{1}{2};\lambda+1,1-\mu,\mu+1;-k^{2}\right)\:,\label{eq:Gegenbauer_2F3}\\
 & = & \sqrt{\pi}\Gamma\left(\lambda+\frac{1}{2}\right)\,_{2}\tilde{F}_{3}\left(\frac{1}{2},\frac{1}{2};\lambda+1,1-\mu,\mu+1;-k^{2}\right)\nonumber 
\end{eqnarray}
where the second form involving the regularized hypergeometric function
\cite{wolfram.com/07.26.26.0001.01} will be needed in numerical evaluations
to avoid infinities and/or indeteminancies whenever $\mu$ is a nonzero
integer. The corresponding Pochhammer symbols in the denominator of
(\ref{eq:H_as_1F2_1F2_Gegen}) will likewise cause numerical problems
unless we switch the $\,_{1}F_{2}$ functions to regularized hypergeometric
functions\cite{wolfram.com/07.22.26.0001.01}. Rewriting those Pochhammer
symbols as ratios of gamma functions,

\begin{equation}
\left(c\right)_{g}=\frac{\Gamma\left(g+c\right)}{\Gamma\left(c\right)}\:,\label{eq:pocc2gamma}
\end{equation}
and simplifying the result gives

\noindent\textbf{Corollary 3a.} For any values of \emph{k}, $\mu$,
and for $\lambda\neq0$ , 

\begin{equation}
\begin{array}{ccc}
&  &
\,_{2}\tilde{F}_{3}\left(\frac{1}{2},\frac{1}{2};\lambda+1,1-\mu,\mu+1;-k^{2}\right) \\
& & \\
& = & {\displaystyle \sum_{L=0}^{\infty}\frac{(-1)^{2L}k^{4L}2^{2L+2\lambda+1}(2L+\lambda)\Gamma\left(L+\frac{1}{2}\right)\Gamma\left(\lambda+\frac{1}{2}\right)\left(\left(\lambda+\frac{1}{2}\right)_{2L}\right){}^{2}\Gamma(2L+\lambda+1)^{2}\Gamma(2L+2\lambda)}{\pi\Gamma(L+1)\Gamma(4L+2\lambda+1)^{2}}}\\
 & \times & \,_{1}\tilde{F}_{2}\left(L+\frac{1}{2};2L+\lambda+1,L-\mu+1;-\frac{k^{2}}{4}\right)\:.\\
 & \times & \,_{1}\tilde{F}_{2}\left(L+\frac{1}{2};2L+\lambda+1,L+\mu+1;-\frac{k^{2}}{4}\right)
\end{array}\label{eq:Corollary_3a}
\end{equation}

When either of $\mu$ or $\nu$ is $-\frac{1}{2}$ the order of the
hypergeometric function on the left-hand side of (\ref{eq:theorem_3})
is also reduced since the second or third lower parameter equals the
first upper parameter $\left(\frac{1}{2}\right)$:

\noindent\textbf{Corollary 3b.} For any values of \emph{k}, $\mu$,
and for $\lambda\neq0$ , 

\begin{equation}
\begin{array}{ccc}
 & & 
\,_{2}F_{3}\left(\frac{\mu}{2}+\frac{1}{4},\frac{\mu}{2}+\frac{3}{4};\lambda+1,\mu+\frac{1}{2},\mu+1;-k^{2}\right)  \\
& & \\
& = & {\displaystyle \frac{2^{\mu-\frac{1}{2}}\Gamma(\lambda+1)\Gamma(\mu+1)}{\Gamma\left(\lambda+\frac{1}{2}\right)\Gamma(\lambda)^{2}}}{\displaystyle \sum_{L=0}^{\infty}\frac{(-1)^{2L}k^{2L-\lambda}\left(\left(\lambda+\frac{1}{2}\right)_{2L}\right){}^{2}2^{6L-\lambda-\mu+\frac{3}{2}}\Gamma(2L+\lambda+1)\Gamma(2L+2\lambda)}{(2L)!(2L+\lambda)\left((2\lambda)_{2L}\right){}^{2}\left((2L+2\lambda)_{2L}\right){}^{2}\left(L+\frac{1}{2}\right)_{\mu+\frac{1}{2}}}}\\
 & \times & J_{2L+\lambda}(k)\,_{1}F_{2}\left(L+\frac{1}{2};2L+\lambda+1,L+\mu+1;-\frac{k^{2}}{4}\right)\:.
\end{array}\:.\label{eq:Corollary_3b}
\end{equation}
An alternative form may be had by replacing the Bessel function by
a $\,_{0}F_{1}$ hypergeometric function using (\ref{eq:J_to_0F1}).

If $\lambda=\frac{\mu}{2}+\frac{\nu}{2}$ the order of the hypergeometric
function on the left-hand side of (\ref{eq:theorem_3}) also is reduced
since the first lower parameter equals the third upper parameter:

\noindent\textbf{Corollary 3ci.} For any values of \emph{k}, $\mu$,
and $\nu$, 

\begin{equation}
\begin{array}{ccc}
& &
\,_{2}F_{3}\left(\frac{1}{2},\frac{\mu}{2}+\frac{\nu}{2}+\frac{1}{2};\mu+1,\nu+1,\mu+\nu+1;-k^{2}\right)  \\
& & \\
& = & {\displaystyle \frac{2^{\mu+\nu}\Gamma(\mu+1)\Gamma(\nu+1)\Gamma\left(\frac{\mu}{2}+\frac{\nu}{2}+1\right)}{\sqrt{\pi}\Gamma\left(\frac{\mu}{2}+\frac{\nu}{2}\right)^{2}\Gamma\left(\frac{\mu}{2}+\frac{\nu}{2}+\frac{1}{2}\right)}}\\
 & \times & {\displaystyle \sum_{L=0}^{\infty}\frac{(-1)^{2L}k^{4L}2^{4L-2\left(\frac{\mu}{2}+\frac{\nu}{2}\right)-\mu-\nu+1}\left(\left(\frac{\mu}{2}+\frac{\nu}{2}+\frac{1}{2}\right)_{2L}\right){}^{2}\Gamma(2L+\mu+\nu)}{(2L)!\left(L+\frac{1}{2}\right)_{\mu+\frac{1}{2}}\left(L+\frac{1}{2}\right)_{\nu+\frac{1}{2}}\left(2L+\frac{\mu}{2}+\frac{\nu}{2}\right)\left((\mu+\nu)_{2L}\right){}^{2}\left((2L+\mu+\nu)_{2L}\right){}^{2}}}\\
 & \times & \,_{1}F_{2}\left(L+\frac{1}{2};L+\mu+1,2L+\frac{\mu}{2}+\frac{\nu}{2}+1;-\frac{k^{2}}{4}\right)\\
 & \times & \,_{1}F_{2}\left(L+\frac{1}{2};2L+\frac{\mu}{2}+\frac{\nu}{2}+1,L+\nu+1;-\frac{k^{2}}{4}\right)
\end{array}\:.\label{eq:Corollary_3ci}
\end{equation}
If either of $\mu$ or $\nu$ is $-\frac{1}{2}$, the order of the
hypergeometric function on the left-hand side of (\ref{eq:Corollary_3ci})
is reduced since the first or second lower parameter equals the first
upper parameter $\left(\frac{1}{2}\right)$:

\noindent\textbf{Corollary 3cii.} For any values of \emph{k} and
$\mu$, 

\begin{equation}
\begin{array}{ccc}
& & 
\,_{1}F_{2}\left(\frac{\mu}{2}+\frac{1}{4};\mu+\frac{1}{2},\mu+1;-k^{2}\right)  \\
& & \\
& = & {\displaystyle \sum_{L=0}^{\infty}\frac{(-1)^{2L}2^{6L-\frac{\mu}{2}+\frac{5}{4}}\Gamma\left(\frac{\mu}{2}+\frac{3}{4}\right)\Gamma(\mu+1)k^{2L-\frac{\mu}{2}+\frac{1}{4}}\left(\left(\frac{\mu}{2}+\frac{1}{4}\right)_{2L}\right){}^{2}\Gamma\left(2L+\frac{\mu}{2}+\frac{3}{4}\right)\Gamma\left(2L+\mu-\frac{1}{2}\right)}{(2L)!\left(2L+\frac{\mu}{2}-\frac{1}{4}\right)\Gamma\left(\frac{\mu}{2}-\frac{1}{4}\right)^{2}\Gamma\left(\frac{\mu}{2}+\frac{1}{4}\right)\left(L+\frac{1}{2}\right)_{\mu+\frac{1}{2}}\left(\left(\mu-\frac{1}{2}\right)_{2L}\right){}^{2}\left(\left(2L+\mu-\frac{1}{2}\right)_{2L}\right){}^{2}}}\\
 & \times & J_{2L+\frac{\mu}{2}-\frac{1}{4}}(k)\,_{1}F_{2}\left(L+\frac{1}{2};2L+\frac{\mu}{2}+\frac{3}{4},L+\mu+1;-\frac{k^{2}}{4}\right)
\end{array}\:.\label{eq:Corollary_3cii}
\end{equation}
An alternative form may be had by replacing the Bessel function by
a $\,_{0}F_{1}$ hypergeometric function using (\ref{eq:J_to_0F1}).

In a similar vein, if $\lambda=\frac{\mu}{2}+\frac{\nu}{2}-\frac{1}{2}$
the order of the hypergeometric function on the left-hand side of
(\ref{eq:theorem_3}) also is reduced since the first lower parameter
equals the second upper parameter:

\noindent\textbf{Corollary 3di.} For any values of \emph{k}, $\mu$,
and $\nu$, 

\begin{equation}
\begin{array}{ccc}
& &
_{2}F_{3}\left(\frac{1}{2},\frac{\mu}{2}+\frac{\nu}{2}+1;\mu+1,\nu+1,\mu+\nu+1;-k^{2}\right) \\
& & \\
 & = & {\displaystyle \frac{\left(2^{2(\mu+\nu)-3}(\mu+\nu-1)\Gamma(\mu+1)\Gamma(\nu+1)\right)}{\pi\Gamma(\mu+\nu-1)}} \\
& & \\
 & \times & {\displaystyle \sum_{L=0}^{\infty}\frac{(-1)^{2L}k^{4L}2^{4L-2\mu-2\nu+2}\left(\left(\frac{\mu}{2}+\frac{\nu}{2}\right)_{2L}\right){}^{2}\Gamma(2L+\mu+\nu-1)}{(2L)!\left(L+\frac{1}{2}\right)_{\mu+\frac{1}{2}}\left(L+\frac{1}{2}\right)_{\nu+\frac{1}{2}}\left(2L+\frac{\mu}{2}+\frac{\nu}{2}-\frac{1}{2}\right)\left((\mu+\nu-1)_{2L}\right){}^{2}\left((2L+\mu+\nu-1)_{2L}\right){}^{2}}}\\
 & \times & \,_{1}F_{2}\left(L+\frac{1}{2};L+\mu+1,2L+\frac{\mu}{2}+\frac{\nu}{2}+\frac{1}{2};-\frac{k^{2}}{4}\right)\\
 & \times & \,_{1}F_{2}\left(L+\frac{1}{2};2L+\frac{\mu}{2}+\frac{\nu}{2}+\frac{1}{2},L+\nu+1;-\frac{k^{2}}{4}\right)
\end{array}\:.\label{eq:Corollary_3di}
\end{equation}

Finally, if either of $\mu$ or $\nu$ is $-\frac{1}{2}$ the order
of the hypergeometric function on the left-hand side of (\ref{eq:Corollary_3di})
is reduced since the first or second lower parameter equals the first
upper parameter $\left(\frac{1}{2}\right)$:

\noindent\textbf{Corollary 3dii.} For any values of \emph{k} and
$\mu$, 

\begin{equation}
\begin{array}{ccc}
& &
\,_{1}F_{2}\left(\frac{\mu}{2}+\frac{3}{4};\mu+\frac{1}{2},\mu+1;-k^{2}\right) \\
& & \\
& = & {\displaystyle \frac{2^{\mu-\frac{1}{2}}\Gamma\left(\frac{\mu}{2}+\frac{1}{4}\right)\Gamma(\mu+1)}{\Gamma\left(\frac{\mu}{2}-\frac{1}{4}\right)\Gamma\left(\frac{\mu}{2}-\frac{3}{4}\right)^{2}}} \\
& & \\
 & \times & {\displaystyle \sum_{L=0}^{\infty}\frac{(-1)^{2L}2^{6L-\frac{3\mu}{2}+\frac{9}{4}}k^{2L-\frac{\mu}{2}+\frac{3}{4}}\left(\left(\frac{\mu}{2}-\frac{1}{4}\right)_{2L}\right){}^{2}\Gamma\left(2L+\frac{\mu}{2}+\frac{1}{4}\right)\Gamma\left(2L+\mu-\frac{3}{2}\right)}{(2L)!\left(2L+\frac{\mu}{2}-\frac{3}{4}\right)\left(L+\frac{1}{2}\right)_{\mu+\frac{1}{2}}\left(\left(\mu-\frac{3}{2}\right)_{2L}\right){}^{2}\left(\left(2L+\mu-\frac{3}{2}\right)_{2L}\right){}^{2}}}\\
 & \times & J_{2L+\frac{\mu}{2}-\frac{3}{4}}(k)\,_{1}F_{2}\left(L+\frac{1}{2};2L+\frac{\mu}{2}+\frac{1}{4},L+\mu+1;-\frac{k^{2}}{4}\right)
\end{array}\:.\label{eq:Corollary_3dii}
\end{equation}
An alternative form may be had by replacing the Bessel function by
a $\,_{0}F_{1}$ hypergeometric function using (\ref{eq:J_to_0F1}).

One might hope for another such reduction of (\ref{eq:theorem_3})
for $\lambda=-\frac{1}{2}$, but the ratio on the right-hand side,
\begin{equation}
\lim_{\lambda\to-\frac{1}{2}}\,\frac{\Gamma(2(L+\lambda))}{\Gamma\left(\lambda+\frac{1}{2}\right)}=\left\{ \begin{array}{ccc}
-\frac{1}{2} &  & L=0\\
0 &  & L>0
\end{array}\right.\:,\label{eq:limit}
\end{equation}
truncates the hoped-for series to one term, the well-known result
\cite{Luke1} (p. 216 No.
6.2.7.39) (\ref{eq:JJ=000020to=0000202F3}). Although the order of
the hypergeometric function will not be reduced, one may continue
on in this fashion by taking $\lambda=-\frac{3}{2}$, and using
\begin{equation}
\lim_{\lambda\to-\frac{3}{2}}\,\frac{\Gamma(2(L+\lambda))}{\Gamma\left(\lambda+\frac{1}{2}\right)}=\left\{ \begin{array}{ccc}
\frac{1}{12} &  & L=0\\
0 &  & L>0
\end{array}\right.\:,\label{eq:limit-1}
\end{equation}
to show that
\begin{eqnarray}
& & \hspace{-1.0cm} \,_{3}F_{4}\left(\frac{1}{2},\frac{\mu}{2}+\frac{\nu}{2}+\frac{1}{2},\frac{\mu}{2}+\frac{\nu}{2}+1;-\frac{1}{2},\mu+1,\nu+1,\mu+\nu+1;-z^{2}\right)  \nonumber  \\
& = & 2^{\mu+\nu}\Gamma(\mu+1)\Gamma(\nu+1)z^{-\mu-\nu}\left(J_{\mu}(z)+zJ_{\mu+1}(z)\right)\left(J_{\nu}(z)+zJ_{\nu+1}(z)\right)\:,\label{eq:lam-3/2}
\end{eqnarray}
which is not among the several simlar relations in Prudnikov, Brychkov,
and Marichev \cite{PBM3} (p. 654 No. 8.4.19..21-2). One may continue in this fashion with negative
half integer values of $\lambda$ of increasing magnitude.
 
.

\section{Conclusions}

We have found that certain $\,_{3}F_{4}$ generalized hypergeometric
functions can be expanded in sums of pair products of $\,_{1}F_{2}$
functions, found by expanding products of Bessel functions in a series
of Chebyshev or Gegenbauer polynomials and using their orthogonality
relations to reduce these double sums to single sums. But that product
my also be directly integrated with no expansion.

In special cases, the $\,_{3}F_{4}$ hypergeometric functions may
be reduced to $\,_{2}F_{3}$ functions, and of those, further special
cases reduce the $\,_{2}F_{3}$ functions to $\,_{1}F_{2}$ functions,
and the sums to products of $\,_{0}F_{1}$ (Bessel) and $\,_{1}F_{2}$
functions. We have, thus, developed hypergeometric function summation
theorems beyond those expressible as pair-products of generalized
Whittaker functions, $\,_{2}F_{1}$ functions, and $\,_{3}F_{2}$
functions into the realm of $\,_{p}F_{q}$ functions where $p<q$
for both the summand and terms in the series.

\end{document}